\documentclass[twoside,a4paper,12pt,reqno]{article} 

\usepackage[utf8]{inputenc} 

\usepackage{pgf,tikz}
\usetikzlibrary{arrows,automata,backgrounds,calendar,chains,matrix,mindmap,patterns,petri,shadows,shapes.geometric,shapes.misc,spy,trees}

\usepackage{geometry} 
\usepackage{graphicx} 

\usepackage{booktabs} 
\usepackage{array} 
\usepackage{paralist} 
\usepackage{verbatim} 
\usepackage{subfig} 
\usepackage{amsfonts} 
\usepackage{amsmath} 
\usepackage{amsthm} 
\usepackage{fancyref} 
\usepackage{hyperref}
\usepackage{pst-eucl}
\usepackage{pstricks}
\usepackage{authblk}
\usepackage{microtype}

\textheight=\dimexpr
  \ifcase 12pt
    \or 
    \or 
    \or 
    \or 
    \or 
    \or 
    \or 
    \or 
    57\or 
    52\or 
    48\or 
    44\or 
    41\fi 
  \baselineskip+\topskip\relax
	
\usepackage{abstract}

\usepackage{fancyhdr}
\usepackage{sectsty}
\allsectionsfont{\sffamily\mdseries\upshape} 

\usepackage[nottoc,notlof,notlot]{tocbibind} 
\usepackage[titles,subfigure]{tocloft} 

\usepackage[square,sort,comma,numbers]{natbib}

\newtheorem{theorem}{Theorem}[section]
\newtheorem*{theorem*}{Theorem}
\newtheorem{corollary}[theorem]{Corollary}
\newtheorem{lemma}[theorem]{Lemma}
\newtheorem*{lemma*}{Lemma}
\newtheorem{proposition}[theorem]{Proposition}
\newtheorem*{proposition*}{Proposition}
\newtheorem{conjecture}[theorem]{Conjecture}
\newtheorem*{conjecture*}{Conjecture}
\newenvironment{claim}[1]{\par\noindent\textbf{Claim:}\space#1}{}

\theoremstyle{definition}
\newtheorem{example}[theorem]{Example}
\newtheorem*{example*}{Example}
\newtheorem{definition}[theorem]{Definition}
\newtheorem*{definition*}{Definition}

\newtheorem*{fact*}{Fact}
\newtheorem{remark}[theorem]{Remark}
\newtheorem*{remark*}{Remark}

\numberwithin{equation}{section}
\newcommand{\R}{\mathbb{R}}
\newcommand{\Sp}{\mathbb{S}}
\newcommand{\Q}{\mathbb{Q}}

\newcommand{\Span}{\mbox{span}}
\newcommand{\noco}{\nabla^\perp}
\newcommand{\lcco}{\nabla}

\newcommand{\virg}{\mbox{, }}
\def\al{\alpha}
\def\delt{\delta}
\def\lamb{\lambda}
\def\tlamb{\tilde{\lamb}}
\def\sigm{\sigma}

\def\Lamb{\Lambda}
\def\vphi{\varphi}
\def\thff{\mbox{\normalfont{III}}}
\def\trac{\mbox{\normalfont{trace\hspace{2pt}}}}

\def \riar{\rightarrow}
\def\<{\langle}
\def\>{\rangle}
\def\Ric{\mbox{Ric}}
\def\U{\mathcal{U}}
\def\beeq{\begin{equation}}
\def\eneq{\end{equation}}
\def\beeqs{\begin{eqnarray*}}
\def\eneqs{\end{eqnarray*}}
\def\besp{\begin{split}}
\def\ensp{\end{split}}
\def\beth{\begin{theorem}}
\def\enth{\end{theorem}}
\def\bepr{\begin{proof}}
\def\enpr{\end{proof}}
\def\beths{\begin{theorem*}}
\def\enths{\end{theorem*}}
\def\becor{\begin{corollary}}
\def\encor{\end{corollary}}
\def\bere{\begin{remark}}
\def\enre{\end{remark}}
\def\beres{\begin{remark*}}
\def\enres{\end{remark*}}
\def\bele{\begin{lemma}}
\def\enle{\end{lemma}}
\def\beles{\begin{lemma*}}
\def\enles{\end{lemma*}}
\def\bepro{\begin{proposition}}
\def\enpro{\end{proposition}}
\def\bepros{\begin{proposition*}}
\def\enpros{\end{proposition*}}
\def\becl{\begin{claim}}
\def\encl{\end{claim}}
\def\beex{\begin{example}}
\def\enex{\end{example}}
\def\beexs{\begin{example*}}
\def\enexs{\end{example*}}
\def\beco{\begin{conjecture}}
\def\enco{\end{conjecture}}
\def\becos{\begin{conjecture*}}
\def\encos{\end{conjecture*}}
\def\bede{\begin{definition}}
\def\ende{\end{definition}}
\def\bedes{\begin{definition*}}
\def\endes{\end{definition*}}

\title{SUBMANIFOLDS WITH HOMOTHETIC GAUSS MAP IN CODIMENSION TWO}
\author{GUILHERME MACHADO DE FREITAS\thanks{Partially supported by CNPq-Brazil.}}
\date{}

\begin{document}
\maketitle

\begin{abstract}
\noindent
ABSTRACT. In this article we derive a complete classification of all submanifolds in space forms with codimension two for which the Gauss map is homothetic.
\end{abstract}

\section{Introduction and statement of the main results}\label{sec:ismr}
Since the very beginning of differential geometry the Gauss map has played an important role in surface theory. A natural generalization of this classical map for an isometric immersion $f:M^n\riar\R^{n+p}$ of an $n$-dimensional Riemannian manifold into the ($n+p$)-dimensional Euclidean space is defined by assigning to every point $x\in M^n$ its tangent space $T_xM$. The Gauss map $\phi:M^n\riar G_n(\R^{n+p})$ into the Grassmannian $G_n(\R^{n+p})$ of $n$-subspaces of $\R^{n+p}$ obtained this way has been extensively studied, and a beautiful survey on results concerning $\phi$ and on alternative definitions of the Gauss map of $f$ can be found in \cite{MR609562}. In this paper we will mainly consider the pullback $\thff$ of the canonical Riemannian metric on $G_n(\R^{n+p})$ (regarded as a symmetric space) via $\phi$, which is called the \emph{third fundamental form} of $f$. It is very natural to pose the following Main Problem:

\vspace{\baselineskip}
\emph{Find all Euclidean submanifolds for which the Gauss map is homothetic (i.e., $\thff$ is a constant multiple of the Riemannian metric on $M^n$).}

\vspace{\baselineskip}
Due to \cite{MR0234388}, $\thff$ can be written in terms of the second fundamental form $\al:TM\times TM\riar N_fM$ of $f$ as
\beeq\label{tffd}
\thff(X,Y)=\sum_{i=1}^{n}\<\al(X,X_i),\al(Y,X_i)\>,
\eneq
where $X_1,\dots,X_n$ is any orthonormal tangent frame. In terms of the mean curvature vector $H$ of $f$ and of the Ricci tensor $\Ric$ of $M^n$, the Gauss equation leads to the invariant description
\beeq\label{tffmcvrt}
\thff(X,Y)=n\<\al(X,Y),H\>-\Ric(X,Y).
\eneq
Notice that for curves, i.e., $n=1$, we have $\thff=\kappa^2\<\cdot,\cdot\>$, where $\kappa$ is the curvature function. Thus, a curve has homothetic Gauss map if and only if it has constant curvature, so that we can assume $n\geq2$.

We obtain from (\ref{tffmcvrt}) a strong connection between our Main Problem and minimal Einstein submanifolds of Euclidean space, namely, a minimal immersion $f:M^n\riar\R^{n+p}$ has homothetic Gauss map if and only if $M^n$ is an Einstein manifold (for $n=2$, by an Einstein surface we mean a surface with constant Gaussian curvature). Another interesting consequence of this equation is the fact that minimal Einstein submanifolds of Euclidean spheres also have homothetic Gauss map. Indeed, minimality in the sphere easily implies that the shape operator in the direction of $H$ is a constant multiple of the identity map. There are important examples of minimal Einstein submanifolds in spheres, the so-called \emph{Veronese embeddings} corresponding to the nonzero eigenvalues of the Laplacian on an irreducible compact symmetric space. They are natural generalizations of the classical Veronese surface $\mathbb{RP}^2\hookrightarrow\Sp^4\subset\R^5$. See \cite{MR0278318}, \cite{MR775143}, \cite{MR0407774} and Chapter 4, §5-6, of \cite{MR749575} for the definition and other concrete examples.

In \cite{MR1066578}, a partial answer to the Main Problem was obtained by Nölker under the assumption of flat normal bundle. Under this restriction, the only non-totally geodesic solutions are Riemannian products of totally umbilical submanifolds with mean curvature vectors of the same constant length, i.e., Euclidean round spheres or curves of constant curvature. Observe that without the assumption of flat normal bundle the Veronese surface provides a counterexample to Nölker's theorem already in codimension three. Nevertheless, we show that such assumption can be dropped in codimension two. Throughout this paper we agree that a \emph{round sphere} $\Sp^n(r)\subset\R^N$ is an $n$-dimensional totally umbilical submanifold of radius $r$, even for $n=1$.

\beth\label{mine}
Let $M^n$ be an $n$-dimensional connected Riemannian manifold, $n\geq2$, and let $f:M^n\riar\R^{n+2}$ be an isometric immersion. Then $\thff=\frac{1}{r^2}\<\cdot,\cdot\>$ with $r>0$ if and only if $f(M^n)$ is (an open subset of) either a round sphere $\Sp^n(r)\subset\R^{n+1}\subset\R^{n+2}$ or a product of two round spheres $\Sp^m(r)\times\Sp^{n-m}(r)\subset\R^{m+1}\times\R^{n-m+1}=\R^{n+2}$.
\enth

As a consequence, we have that there is no substantial irreducible codimension two Euclidean submanifold with homothetic Gauss map (except curves of constant curvature in $\R^3$).

The key fact that the third fundamental form can be written in terms of the second fundamental form allows us to naturally extend our Main Problem for isometric immersions into real space forms $\Q_c^N$ of nonzero curvature.

A version of Nölker's theorem for the case $c\neq0$ can be easily obtained, based on the notion of extrinsic products of isometric immersions; cf. Remark in Section 1 of \cite{MR1066578}. Let us recall this construction.

Let us regard the space form $\Q_c^{N-1}$ as
\beeqs
\Q_c^{N-1}=\big\{X=(x_1,\dots,x_N)\in\mathbb{E}^N:\<X,X\>=\frac{1}{c}\big\}.
\eneqs
where $\mathbb{E}^N$ denotes either the Euclidean space $\R^N$ or the Lorentz space $\mathbb{L}^N$ according to whether $c>0$ or $c<0$, respectively, and $x_1>0$ in the latter case. Given an orthogonal decomposition
\beeqs
\mathbb{E}^N=\mathbb{E}^{m_0}\times\R^{m_1}\times\cdots\times\R^{m_k}
\eneqs
and immersions $f_0:M_0^{n_0}\riar\Q_{c_0}^{m_0-1}\subset\mathbb{E}^{m_0}$ and $f_i:M_i^{n_i}\riar\Sp_{c_i}^{m_i-1}\subset\R^{m_i}\virg 1\leq i\leq k$, then the image $f(M^n)$ of their product immersion $f=f_0\times f_1\times\cdots\times f_k:M^n=M_0^{n_0}\times M_1^{n_1}\times\cdots\times M_k^{n_k}\riar\mathbb{E}^N=\mathbb{E}^{m_0}\times\R^{m_1}\times\cdots\times\R^{m_k}$ given by
\beeqs
f(p_1,\dots,p_k)=(f_0(p_0),f_1(p_1),\dots,f_k(p_k))
\eneqs
is contained in the space form $\Q_c^{N-1}$ of constant sectional curvature $c=\big(\sum_{i=0}^{k}\frac{1}{c_i}\big)^{-1}$, provided that $\sum_{i=0}^{k}\frac{1}{c_i}\neq0$. If $f$ is regarded as an immersion into $\Q_c^{N-1}$, then it is called the \emph{extrinsic product} of $f_0\virg f_1\virg\dots\virg f_k$. On the other hand, consider now an orthogonal decomposition $\R^{N-1}=\R^{m_1}\times\cdots\times\R^{m_k}$ and isometric immersions $f_i:M_i\riar\R^{m_i}\virg 1\leq i\leq k$. If we regard each $f_i\virg 1\leq i\leq k$, as an isometric immersion into $\R^{N-1}$ and consider its composition $\tilde{f_i}=j\circ f_i$ with the umbilical inclusion $j:\R^{N-1}\riar\mathbb{H}_c^N$, then we also say that $\tilde{f}=j\circ(f_1\times\cdots\times f_k)$ is the extrinsic product of $\tilde{f}_1\virg\dots\virg\tilde{f}_k$.

Under the above notation, Nölker's argument can be generalized in space forms to show that every non-totally geodesic isometric immersion into space forms with flat normal bundle and homothetic Gauss map is an extrinsic product of either totally umbilical isometric immersions $f_0\virg f_1\virg\dots\virg f_k$ or $\tilde{f_1}\virg\dots\virg\tilde{f_k}$, where the mean curvature vectors $H_i$ of $f_i$ have all the same constant length in the latter case, $1\leq i\leq k$, and $\left\|H_i\right\|=\sqrt{\rho^2-c_i}$ for some $\rho>0$ in the former, $0\leq i\leq k$.

We say that an isometric immersion in space forms is \emph{irreducible} if it does not split as an extrinsic product as above. Our next result shows that the only substantial irreducible solution of the Main Problem for codimension two submanifolds in space forms is the Veronese surface in the 4-sphere. This together with the preceding discussion provides a complete classification of all submanifolds in space forms with codimension two and homothetic Gauss map.

\beth\label{mine2}
Let $M^n$ be an $n$-dimensional connected Riemannian manifold, $n\geq2$, and let $f:M^n\riar\Q_c^{n+2}$ be a substantial irreducible isometric immersion with homothetic Gauss map. Then $n=2\virg c>0$ and $f(M^n)$ is (an open subset of) the Veronese surface $\mathbb{RP}^2\hookrightarrow\Sp_c^4$.
\enth

\vspace{\baselineskip}
\emph{Acknowledgments}. The author is grateful to his Ph.D advisor at IMPA, Prof. Luis A. Florit, for his constant advice and encouragement. The author also wishes to thank Profs. Marcos Dajczer, Ruy Tojeiro and Antonio Di Scala for helpful discussions and comments.

\section{Minimal Einstein submanifolds}
Here, we state some results related to minimal Einstein submanifolds which will be necessary for the proofs of Theorems \ref{mine} and \ref{mine2}.

As previously mentioned, we have the following fact.

\bepro\label{emes}
Let $f:M^n\riar\R^{n+p}$ be a minimal immersion, with $n\geq2$. Then, $f$ has homothetic Gauss map if and only if $M^n$ is an Einstein manifold.
\enpro

By combining a result of Osserman-Chern \cite{MR0226514} for the Gauss map and a result of Calabi \cite{MR0057000} for Riemann surfaces in the complex projective spaces, we know that the hyperbolic plane can not be minimally immersed into a Euclidean space even locally. In other words, we have

\beth\label{mihe}
Every minimal surface in Euclidean space with constant Gauss curvature must be totally geodesic.
\enth

The next conjecture, due to Di Scala \cite{MR2000031}, is the higher-dimensional version of
the previous result.

\beco\label{disc}
Let $M^n$ be an Einstein manifold, with $n\geq3$. Then, any minimal isometric immersion $f:M^n\riar\R^{n+p}$ must be totally geodesic.
\enco

According to the main result of \cite{MR2000031}, the conjecture is true if $M^n$ is also Kähler. Furthermore, under the assumption of flat normal bundle, it follows as a corollary
of Nölker's theorem and Proposition \ref{emes}. In \cite{MR1255769}, Matsuyama presented a general proof in codimension two. His result is stated below.
\beth\label{mats}
Let $M^n$ be an Einstein manifold, with $n\geq3$. Then, any minimal isometric immersion $f:M^n\riar\R^{n+2}$ with codimension two must be totally geodesic.
\enth
 Theorems \ref{mihe} and \ref{mats} are also true for minimal Einstein submanifolds of hyperbolic space (see \cite{MR787964} and \cite{MR1255769}).

In the sphere, though, the situation is different. In \cite{MR706526}, Kenmotsu has provided a complete classification of the minimal surfaces with constant Gaussian curvature in the 4-sphere. The only non-totally geodesic ones are the Clifford torus and the Veronese surface. Notice that the Clifford torus is reducible in the sense of extrinsic products. In higher dimension, Matsuyama \cite{MR1255769} classified the minimal Einstein submanifolds with codimension two in the sphere. The only such submanifolds are products of up to three spheres of the same dimension and radius.

\section{The algebraic decomposition}

In this section, we prove some algebraic results that will play a key role in the proofs of Theorems \ref{mine} and \ref{mine2}.

We write $I_V$ for the identity automorphism on a vector space $V$. The spectrum of a self-adjoint operator $A$ and the eigenspace associated to the eigenvalue $\lamb$ are denoted by $\Lamb_A$ and $E_A(\lamb)$, respectively. For convenience, we set $E_A(\lamb)=\{0\}$ for $\lamb\in\R\setminus\Lamb_A$.

Let $V$ and $W$ be real vector spaces of finite dimension with positive definite inner products and let $\al:V\times V\riar W$ be a symmetric bilinear form. For any given $\xi\in W$, we define the \emph{shape operator} $A_\xi:V\riar V$ of $\al$ with respect to $\xi$ by
\beeqs
\<A_\xi X,Y\>=\<\al(X,Y),\xi\>.
\eneqs
We say that $\al$ is \emph{adapted} to an orthogonal decomposition $V=E_1\oplus\ldots\oplus E_k$ if the subspaces $E_i\virg1\leq i\leq k$, are preserved by all shape operators. Equivalently,
\beeqs
\al(E_i,E_j)=0\virg\forall1\leq i\neq j\leq k.
\eneqs 

Finally, a bilinear form $\vphi:V\times V\riar W$ is said to be \emph{umbilical} if there exists a vector $\xi\in W$ such that
\beeqs
\vphi(X,Y)=\<X,Y\>\xi
\eneqs
for all $X\virg Y\in V$. We start with a useful criterion for umbilical bilinear forms. 

\bele\label{alle}
Let $V$ and $W$ be real vector spaces of finite dimension, where $V$ has a positive definite inner product, and let $\vphi:V\times V\riar W$ be a bilinear form such that $\vphi(X,Y)=0$ for all pair of orthonormal vectors $X\virg Y\in V$. Then $\vphi$ is umbilical.
\enle

\bepr
Let $\{X_1,\dots,X_n\}$ be an orthonormal basis of $V$ and set
\beeqs
\vphi_{ij}=\vphi(X_i,X_j)\virg1\leq i\virg j\leq n.
\eneqs
By linearity, all we have to prove is that $\vphi_{ij}=\delt_{ij}\xi$ for some $\xi\in W$. For $i\neq j$, it holds by assumption. For $i=j$, take the orthonormal vectors $X=\frac{1}{\sqrt{2}}(X_i+X_k)\virg Y=\frac{1}{\sqrt{2}}(X_i-X_k)\virg i\neq k$, and use the assumption to conclude that
\beeqs
0=\vphi(X,Y)=\frac{1}{2}(\vphi_{ii}-\vphi_{kk}).
\eneqs
Thus $\xi=\vphi_{ii}$ does not depend on $i$ and our lemma is proved.
\enpr

Next, we study some algebraic implications of having umbilical third fundamental form. We use equation (\ref{tffd}) as an abstract definition of the \emph{third fundamental form} associated to a symmetric bilinear form $\al:V\times V\riar W$.

\bele \label{aldec}
Let $V^n$ and $W^2$ be real vector spaces of dimensions $n$
and 2, respectively, endowed with positive definite inner products, and let $\al:V^n\times V^n\riar W^2$ be a symmetric bilinear form. If the third fundamental form $\thff$ associated to $\al$ is umbilical, then there is an integer $k\virg1\leq k\leq n$, pairwise distinct nonnegative functions $\lamb_j:W^2\riar\R_{\geq0}\virg1\leq j\leq k$, and an orthogonal decomposition
\beeq \label{alde}
V^n=E_1\oplus\ldots\oplus E_k
\eneq
to which $\al$ is adapted and such that the shape operators satisfy
\beeqs
A_\xi^2|_{E_j}=\lamb_j^2(\xi)I_{E_j}
\eneqs
for all $\xi\in W^2$. Moreover, the integer $k$, the functions $\lamb_j\virg1\leq j\leq k$, and the above decomposition are unique (up to permutations).
\enle
\bepr
Observe that $\Lamb_{A_\xi^2}=\{\kappa^2:\kappa\in\Lamb_{A_\xi}\}$. Moreover,
\beeq \label{A2}
E_{A_\xi^2}(\kappa^2)=E_{A_\xi}(\kappa)\oplus E_{A_\xi}(-\kappa)
\eneq
for any $\xi\in W^2$. In particular, $A_\xi$ leaves the eigenspaces $E_{A_\xi^2}(\kappa^2)$ of $A_\xi^2$ invariant.

We can assume that $\thff\neq0$, since, otherwise, $\al=0$ by (\ref{tffd}) and there is nothing to prove. The assumption that $\thff$ is umbilical, say $\thff=\frac{1}{r^2}\<\cdot,\cdot\>$, is equivalent to $A_{\xi_1}^2+A_{\xi_2}^2=\frac{1}{r^2}I_V$, where $\{\xi_1,\xi_2\}$ is any orthonormal basis of $W^2$. In other words, the spectra and eigenspaces of $A_{\xi_1}^2\mbox{ and }A_{\xi_2}^2$ are related by 
\begin{gather*}
\Lamb_{A_{\xi_2}^2}=\big\{\frac{1}{r^2}-\lamb^2:\lamb^2\in\Lamb_{A_{\xi_1}^2}\big\},\\
E_{A_{\xi_2}^2}\big(\frac{1}{r^2}-\lamb^2\big)=E_{A_{\xi_1}^2}(\lamb^2).
\end{gather*}
Thus, both $A_{\xi_1}$ and $A_{\xi_2}$ must leave the eigenspaces of $A_{\xi_i}^2$ invariant, $1\leq i\leq2$. As we are in codimension two, it follows that $\al$ is adapted to the eigendecomposition of $A_{\xi_i}^2$. But since the orthonormal basis $\{\xi_1,\xi_2\}$ of $W^2$ was taken arbitrarily, we conclude that $\al$ is indeed adapted to the eigendecomposition of any $A_\xi^2\virg\xi\in W^2$. In other words,
\beeqs
A_\eta E_{A_\xi^2}(\lamb^2)\subseteq E_{A_\xi^2}(\lamb^2)\virg\forall\xi\virg\eta\in W^2,
\eneqs
where $\lamb^2\in\Lamb_{A_\xi^2}$. In particular, the eigenspaces of each $A_\xi^2$ are invariant under any other $A_\eta^2$. Since the endomorphisms $A_\xi^2$ are self-adjoint, this is equivalent to the existence of a common eigenbasis for the family $\{A_\xi^2:\xi\in W^2\}$. It is now straightforward to verify that the components $E_j$ in our decomposition (\ref{alde}) must be precisely the \emph{eigenspaces of the family} $\{A_\xi^2:\xi\in W^2\}$, i.e., the maximal subspaces of common eigenvectors of all $A_\xi^2\virg\xi\in W^2$. Equivalently,
\beeqs
E_j=\cap_{\xi\in W^2}E_{A_\xi^2}(\lamb_j^2(\xi)),
\eneqs
where the eigenvalues $\lamb_j^2(\xi)\in\Lamb_{A_\xi^2}$ are such that the subspace on the right-hand side of the above equality is nonzero. This concludes the proof of the lemma.
\enpr

\bere \label{nola}
Notice that, if $\thff=\frac{1}{r^2}\<\cdot,\cdot\>$, then $\lamb_j(\xi_1)^2+\lamb_j(\xi_2)^2=\frac{1}{r^2}$ for $1\leq j\leq k$ and every orthonormal basis $\{\xi_1,\xi_2\}$ of $W^2$.
\enre

The idea now is to understand the algebraic structure of $\al$ restricted to each block of decomposition (\ref{alde}).

\bele \label{aebb}
Let $E^m$ and $W^2$ be real vector spaces of dimensions $m$
and 2, respectively, endowed with positive definite inner products, and let $\al:E^m\times E^m\riar W^2$ be a symmetric bilinear form. If there exists a positive function $\lamb:W^2\setminus\{0\}\riar\R_{>0}$ such that the shape operators of $\al$ satisfy
\beeq\label{sbso}
A_\xi^2=\lamb(\xi)^2I_E
\eneq
for every $\xi\in W^2\setminus\{0\}$, then both $\lamb(\xi)\mbox{ and }-\lamb(\xi)$ are eigenvalues of $A_\xi$ and have the same multiplicity (in particular, $m$ is even and $\trac A_\xi=0$).

Furthermore, for any orthonormal basis $\{\xi_1,\xi_2\}$ of $W^2$,  there exist $\rho\virg\sigm\in\R\virg\sigm\geq0$, satisfying $\rho^2+\sigm^2=\lamb(\xi_2)^2$ and a linear map $A:E^+\riar E^-$ such that
\beeq\label{AA*s}
A^*A=\sigm^2I_{E^+}\hspace{20pt}AA^*=\sigm^2I_{E^-}
\eneq
and
\beeq\label{Axi2}
\begin{gathered}
  A_{\xi_1}=\lamb(\xi_1)(\pi^+-\pi^-),\\
	A_{\xi_2}=(\rho I_E+A)\pi^++(-\rho I_E+A^*)\pi^-,
\end{gathered}
\eneq
where $E^\pm=E_{A_{\xi_1}}(\pm\lamb(\xi_1))$ and $\pi^\pm$ is the orthogonal projection $\pi^\pm:E^m\riar E^\pm$.
\enle
\bepr
Take any orthonormal basis $\{\xi_1,\xi_2\}$ of $W^2$. We have that
\beeqs
A_{\xi_1+\xi_2}^2=(A_{\xi_1}+A_{\xi_2})^2=A_{\xi_1}^2+A_{\xi_2}^2+A_{\xi_1} A_{\xi_2}+A_{\xi_2} A_{\xi_1}.
\eneqs
By the assumption, we obtain
\beeq \label{pcmi}
A_{\xi_1} A_{\xi_2}+A_{\xi_2} A_{\xi_1}=\beta I_E
\eneq
for some real number $\beta$.

For simplicity of notation set $\lamb(\xi_1)=\tlamb$. Since $A_{\xi_1}^2=\tlamb^2 I_E$, it follows that $\Lamb_{A_{\xi_1}}\subseteq\{-\tlamb,\tlamb\}$. Write
\beeqs
A_{\xi_1}=\tlamb(\pi^+-\pi^-),\hspace{10pt}A_{\xi_2}=(A^++A)\pi^++(A^-+B)\pi^-
\eneqs
according to the eigendecomposition of $A_{\xi_1}$, where $A^\pm:E^\pm\riar E^\pm$, $A:E^+\riar E^-$ and $B:E^-\riar E^+$. As $A_{\xi_2}$ is self-adjoint, we must have
\beeq \label{A2sa}
B=A^*.
\eneq
Then $A_{\xi_1} A_{\xi_2}+A_{\xi_2} A_{\xi_1}=2\tlamb(A^+\pi^+-A^-\pi^-)$ and (\ref{pcmi}) yields
\beeq \label{A+A-}
A^\pm=\pm\rho I_{E^\pm}
\eneq
with $\rho=\frac{\beta}{2\tlamb}$. Now, it follows using (\ref{A2sa}) and (\ref{A+A-}) that $A_{\xi_2}^2=\rho^2 I_E+A^*A\pi^++AA^*\pi^-$. Since $A^*A$ and $AA^*$ are positive operators, we conclude by using the assumption on the shape operators again that $\lamb(\xi_2)^2-\rho^2\geq0$ and obtain (\ref{AA*s}) for $\sigm^2=\lamb(\xi_2)^2-\rho^2$. If $\sigm=0$, then $A=0$ and thus $A_{\xi_2}=\rho(\pi^+-\pi^-)$. But this implies that $A_\xi=0$ for some $\xi\neq0$, which contradicts the positivity of $\lamb$. Therefore, $\sigm\neq0$ and $A:E^+\riar E^-$ is an isomorphism. In particular, both $\tlamb$ and $-\tlamb$ are eigenvalues of $A_{\xi_1}$ and have the same multiplicity. Since $\xi_1\in W^2\setminus\{0\}$ was taken arbitrarily, the proof is complete. 
\enpr

\bere \label{spfr}
Observe that $\lamb$ is positive if and only if $\mbox{Im }\al=\Span\{\al(X,Y):X\virg Y\in E^m\}$ is two-dimensional. In the case where there is a nonzero vector $\xi\in W^2$ such that $\lamb(\xi)=0$, (\ref{AA*s}) and (\ref{Axi2}) still hold provided that $\xi_1$ is not collinear to $\xi$. However, $A=0$ and then $\dim E^+\neq\dim E^-$ in general.
\enre

We finally compile the information contained in (\ref{AA*s}) and (\ref{Axi2}) by means of certain umbilical bilinear forms derived from $\al$ and $A$. Define $\al_A:E^+\times E^+\riar W^2$ and $\al_{A^*}:E^-\times E^-\riar W^2$ by $\al_A(X,Y)=\al(X,AY)$ and $\al_{A^*}(X,Y)=\al(X,A^*Y)$.
\bele\label{alum}
Let $\al$ and $A$ be as in Lemma \ref{aebb}. Then the bilinear forms $\al|_{E^+\times E^+}$, $\al|_{E^-\times E^-}$, $\al_A$, $\al_{A^*}$ are all umbilical. More precisely, we have
\beeq\label{alumeq}
\begin{aligned}
\al|_{E^+\times E^+}(X,Y)&=\<X,Y\>(\lamb(\xi_1)\xi_1+\rho\xi_2), & \al_A(X,Y)&=\<X,Y\>\sigm^2\xi_2,\\
\al|_{E^-\times E^-}(X,Y)&=-\<X,Y\>(\lamb(\xi_1)\xi_1+\rho\xi_2), & \al_{A^*}(X,Y)&=\<X,Y\>\sigm^2\xi_2
\end{aligned}
\eneq
for all $X\virg Y$ in the corresponding domains.
\enle

\bepr
We argue for $\al_A$, the other cases being similar. Since $E^+$ and $E^-$ are eigenspaces of $A_{\xi_1}$ associated to distinct eigenvalues, it follows that
\beeq\label{alA1}
\<\al_A(X,Y),\xi_1\>=0
\eneq
for any $X\virg Y\in E^+$. On the other hand, (\ref{AA*s}) and (\ref{Axi2}) imply that
\beeqs
\<\al_A(X,Y),\xi_2\>=\<AX,AY\>=\<X,Y\>\sigm^2.
\eneqs
Therefore $\al_A(X,Y)=\<X,Y\>\sigm^2\xi_2$ and the lemma is proved.
\enpr

\section{Proofs of Theorems \ref{mine} and \ref{mine2}} \label{sec:proo}
The main goal of this section is to prove Theorems \ref{mine} and \ref{mine2}. Additionally, we conclude posing a conjecture suggesting a possible complete solution to our Main Problem in arbitrary codimension.

Throughout this section, we denote by $\lcco$ and $\noco$ the Levi-Civita connection of $M^n$ and the normal connection of $f$, respectively.

The following result is of independent interest.
\bepro\label{roco}
Let $f:M^n\riar\R^{n+p}$ be an isometric immersion, and suppose that there exists a totally geodesic submanifold $L$ in $M^n$ such that $\al$ is adapted to $(TL,TL^\perp\cap TM)$. Then $f|_L$ admits a reduction of codimension to $p$.
\enpro

\bepr
Let $(\al,\noco)$ and $(\al_L,{}^L\noco)$ denote the second fundamental forms and normal connections of $f$ and $f|_L$, respectively. In terms of the second fundamental forms, the assumption that $L$ is a totally geodesic submanifold of $M^n$ means that
\beeq\label{sffg}
\al_L=\al|_{TL\times TL}.
\eneq
In particular, we have $N_1L\subseteq N_1M$, where $N_1M$ and $N_1L$ are the first normal spaces of $f$ and $f|_L$, respectively.

The assumption that $\al$ is adapted to $(TL,TL^\perp\cap TM)$ implies that $A_\xi X\in TL$ for all $\xi\in N_fM$ whenever $X\in TL$. In other words,
\beeq\label{sofl}
A_\xi^L=A_\xi|_{TL}
\eneq
for all $\xi\in N_fM$, where $A_\xi^L$ denotes the shape operator of $f|_L$ with respect to $\xi$. Hence, comparing the Weingarten formulas of $f$ and $f|_L$ we obtain
\beeqs
{}^L\noco_{X}\xi=\noco_{X}\xi\virg\forall X\in TL\mbox{ and }\xi\in N_f M.
\eneqs
Therefore, $N_f M$ is a parallel subbundle of rank $p$ of the normal bundle $N_{f|_L}L$ containing $N_1L$. The statement then follows from a well-known fact about reduction of codimension (cf. \cite{MR1075013}).
\enpr

\bere\label{xxxx} Let $\thff_L:TL\times TL\riar\R$ denote the third fundamental form of $f|_L$. Since $\al$ is adapted to $(TL,TL^\perp\cap TM)$, it follows immediately from (\ref{sffg}) that
\beeqs
\thff_L=\thff|_{TL\times TL},
\eneqs
where $\thff$ is the third fundamental form of $f$.
\enre

Let us start to carry out the proofs of Theorems \ref{mine} and \ref{mine2}.
By Lemma \ref{aldec} we get at each $x\in M^n$ an integer $k\virg1\leq k\leq n$, pairwise distinct nonnegative functions $\lamb_j:N_{f}M(x)\riar\R_{\geq0}\virg1\leq j\leq k$, and an orthogonal decomposition
\beeq \label{alded}
T_xM=E_1\oplus\ldots\oplus E_k
\eneq
to which the second fundamental form $\al_x:T_xM\times T_xM\riar N_fM(x)$ is adapted and such that the shape operators satisfy
\beeqs
A_\xi^2|_{E_j}=\lamb_j(\xi)^2I_{E_j}\virg\forall\xi\in N_fM(x).
\eneqs
In particular, since the integer $k$, the functions $\lamb_j\virg1\leq j\leq k$, and the above decomposition are unique up to permutations, we can choose them to be smooth along an open dense subset $U$ of $M^n$. In fact, we first claim that, at each point $x\in M^n$, there exists a normal vector $\xi_x\in N_fM(x)$ such that the numbers $\lamb_1(\xi_x)\virg\dots\virg\lamb_k(\xi_x)$ are pairwise distinct. Suppose otherwise and let $l<k$ be the maximum number such that $\lamb_1(\xi_0)\virg\dots\virg\lamb_l(\xi_0)$ are pairwise distinct for some $\xi_0\in N_fM(x)$. Pick vectors $\xi_i\in N_fM(x)$ for which $\lamb_{l+1}(\xi_i)\neq\lamb_i(\xi_i)$, $1\leq i\leq l$. This is possible, since the functions $\lamb_1\virg\dots\lamb_k$ are pairwise distinct. Now, set $\xi_x=\xi_0+\sum_{i=1}^{l}t_i\xi_i$ and observe that, for a unit vector $X_m\in E_m$, we have
\small
\beeq
\begin{gathered}
\lamb_m^2(\xi_x)=\<A_{\xi_0+\sum\nolimits_{i=1}^{l}t_i\xi_i}^2X_m,X_m\>\\
=\lamb_m^2(\xi_0)+2\sum\nolimits_{i=1}^{l}(\<A_{\xi_0}X_m,A_{\xi_i}X_m\>+
\sum\nolimits_{j\neq i}\<A_{\xi_i}X_m,A_{\xi_j}X_m\>t_j)t_i+\sum\nolimits_{i=1}^{l}\lamb_m^2(\xi_i)t_i^2.
\end{gathered}
\eneq
\normalsize
Thus, $p_i=\lamb_{l+1}^2(\xi_x)-\lamb_i^2(\xi_x)$ is a quadratic polynomial in the variable $t_i$, and hence has at most two zeros, $1\leq i\leq l$. In particular, $p_i\neq0$ for $t_i$ sufficiently small (once the remaining variables $t_j$ with $j\neq i$ have been fixed). Therefore, we can suitably choose $t_1\virg\dots\virg t_l$ such that $\lamb_1(\xi_x)\virg\dots\virg\lamb_l(\xi_x)$ remain pairwise distinct and $p_i\neq0$ for all $i=1\virg\dots\virg l$. But this implies that $\lamb_1(\xi_x)\virg\dots\virg\lamb_l(\xi_x)\virg\lamb_{l+1}(\xi_x)$ are still pairwise distinct, contradicting the maximality of $l$. This concludes the proof of our claim.

Now, extend $\xi_x$ to a smooth unit normal vector field $\xi$ in a neighborhood of $x$. As the number of eigenvalues of $A_{\xi}$ is a lower semi-continuous function, so is $k(x)$. In particular, $k$ is constant along the connected components of an open dense subset $U$ of $M^n$. Furthermore, since $\lamb_1\virg\dots\virg\lamb_k$ and $E_1\virg\dots\virg E_k$ are, respectively, the eigenvalues and eigenspaces of a shape operator by the above argument, we conclude that they are smooth along each connected component of $U$, as we wished.

\vspace{\baselineskip}
To prove Theorem \ref{mine}, it suffices to show that $f$ has flat normal bundle. Suppose otherwise and take a point $x\in U$ at which this property fails. By the Ricci equation, it means that the shape operators $\{A_\xi:\xi\in N_fM(x)\}$ are not simultaneously diagonalizable. Thus, there is at least one index $j\virg1\leq j\leq k$, such that the family $\{A_\xi|_{E_j}:\xi\in N_fM(x)\}$ is not simultaneously diagonalizable. In particular, since we are in codimension two, no $A_\xi|_{E_j}$ with $\xi\neq0$ can vanish identically, so that $\lamb_j(\xi)\neq0$ for every $\xi\in N_fM(x)\setminus\{0\}$ and Lemma \ref{aebb} applies to $\al|_{E_j\times E_j}$. This clearly remains valid in a small neighborhood $U'\subset U$ of $x$.

\bele\label{Ejtg}
$E_j$ is a totally geodesic (hence integrable) distribution on $U'$.
\enle

\bepr
Let $E_i$ be another distribution in decomposition (\ref{alded}). Since $E_j$ and $E_i$ are orthogonal, we can define a tensor $\vphi:E_j\times E_j\riar E_i$ by projecting $\lcco_XY$ orthogonally onto $E_i$, i.e.,
\beeq \label{toge}
\vphi(X,Y)=(\lcco_XY)_{E_i}.
\eneq
All we have to prove is that $\vphi$ vanishes identically. Consider $\U_0\subset U'$ the set where there is a nonzero normal vector $\xi\in N_fM$ such that $\lamb_i(\xi)=0$.

At each point in $U'$, the functions $\lamb_j$ and $\lamb_i$ are distinct, so that we can take a (local) smooth unit normal vector field $\xi_1$ for which $\lamb_j(\xi_1)\neq\lamb_i(\xi_1)$ everywhere. Furthermore, when working in $\U_0$, we choose $\xi_1$ such that $\lamb_i(\xi_1)\neq0$. Let $\{\xi_1,\xi_2\}$ be a smooth orthonormal normal frame. We write $\lamb_j(\xi_1)=\tlamb_j\virg\lamb_i(\xi_1)=\tlamb_i$ for simplicity and denote by $(\rho_j,\sigm_j,A_j:E_j^+\riar E_j^-)$ and $(\rho_i,\sigm_i,A_i:E_i^+\riar E_i^-)$ the triples given by Lemma \ref{aebb} and Remark \ref{spfr} applied to $\al|_{E_j\times E_j}$ and $\al|_{E_i\times E_i}$, respectively (recall that $\sigm_j\neq0$ and $A_j$ is an isomorphism).

In what follows, the fact that $\al$ is adapted to (\ref{alded}), together with (\ref{AA*s}), (\ref{Axi2}) and (\ref{alumeq}), is often used without explicit mention.

Let us define tensors $\vphi_{A_j^\pm}:E_j^\pm\times E_j^\pm\riar E_i$ by
\beeqs
\vphi_{A_j^\pm}(X,Y)&=\vphi(X,A_j^\pm Y),
\eneqs
where we write $A_j^+=A_j$ and $A_j^-=A_j^*$. Since $A_j^\pm:E_j^\pm\riar E_j^\mp$ is an isomorphism, it suffices to show that $\vphi|_{E_j^\pm\times E_j^\pm}$ and $\vphi_{A_j^\pm}$ vanish identically to conclude the proof of the lemma. Our first goal is to show that, since $\al|_{E_j^\pm\times E_j^\pm}$ and $\al_{A_j^\pm}$ are umbilical bilinear forms by Lemma \ref{alum}, the same property holds for $\vphi|_{E_j^\pm\times E_j^\pm}$ and $\vphi_{A_j^\pm}$. The symbol $\mp$ is used when $\pm$ has already appeared in the same context, to indicate the sign opposite to the one represented by the latter.

Using that $\al$ is adapted to (\ref{alded}) together with (\ref{alumeq}), the Codazzi equation for $(Z\in E_i,X\in E_j^\pm,Y\in E_j^\pm:Y\perp X)$ becomes
\beeqs
\al(\lcco_{Z}X,Y)+\al(X,\lcco_{Z}Y)=\al(\lcco_{X}Z,Y)+\al(Z,\lcco_{X}Y).
\eneqs
Taking the inner product of the equation above with $\xi_1$, the pairwise orthogonality of $X\virg Y\virg Z$ yields (recall that $E_j^\pm$ is the eigenspace of $A_{\xi_1}|_{E_j}$ associated to $\pm\tlamb_j$)
\beeqs
\<\vphi(X,Y),(A_{\xi_1}\mp\tlamb_j I_{E_i})Z\>=0.
\eneqs
Since $\pm\tlamb_j\notin\Lamb_{A_{\xi_1}|_{E_i}}$ (after all, $\tlamb_j\neq\tlamb_i$), we have that $A_{\xi_1}|_{E_i}\mp\tlamb_j I_{E_i}$ is an isomorphism of $E_i$ and thus $\vphi(X,Y)=0$ for all orthonormal pair $X\virg Y\in E_j^\pm$. Therefore, it follows from Lemma \ref{alle} that the bilinear form $\vphi|_{E_j^\pm\times E_j^\pm}$ is umbilical. In other words, there exists a vector field $P^\pm\in E_i$ such that
\beeqs
\vphi|_{E_j^\pm\times E_j^\pm}=\<\cdot,\cdot\>P^\pm.
\eneqs
Taking now the inner product of the same equation with $\xi_2$, (\ref{Axi2}) and again the pairwise orthogonality of $X\virg Y\virg Z$ give (we use the above to conclude that the term $\<\vphi(X,Y),(A_{\xi_2}\mp\rho_j I_{E_i})Z\>$ vanishes)
\beeq\label{codulleta}
\<\vphi_{A_j^\pm}(X,Y),Z\>=-\<\lcco_{Z}X,A_j^\pm Y\>-\<\lcco_{Z}Y,A_j^\pm X\>.
\eneq
In particular, as the right-hand side is symmetric in $X\virg Y$, so is the bilinear form $\vphi_{A_j^\pm}$. 

Using that $\al(X,A_j^\pm Y)=\al_{A_j^\pm}(X,Y)=0$ by Lemma \ref{alum}, the Codazzi equation for $(Z\in E_i,X\in E_j^\pm,A_j^\pm Y:Y\in E_j^\pm;Y\perp X)$ yields
\beeqs
\al(\lcco_{Z}X,A_j^\pm Y)+\al(X,\lcco_{Z}A_j^\pm Y)=\al(\lcco_{X}Z,A_j^\pm Y)+\al(Z,\lcco_{X}A_j^\pm Y).
\eneqs
Taking the inner product of the above equation with $\xi_1$ and taking into account that $A_j^\pm Y\in E_j^\mp$, we obtain
\beeq\label{codulAlxi}
\<\vphi_{A_j^\pm}(X,Y),(A_{\xi_1}\pm\tlamb_j I_{E_i})Z\>=\mp2\tlamb_j\<\lcco_ZX,A_j^\pm Y\>,
\eneq
which alongside the symmetry of $\vphi_{A_j^\pm}$ gives $\<\lcco_{Z}X,A_j^\pm Y\>=\<\lcco_{Z}Y,A_j^\pm X\>$. This and (\ref{codulleta}) then yield
\beeq\label{codulAlss}
\<\vphi_{A_j^\pm}(X,Y),Z\>=-2\<\lcco_{Z}X,A_j^\pm Y\>.
\eneq
Now, multiply (\ref{codulAlss}) by $\mp\tlamb_j$ and add the result to (\ref{codulAlxi}), to obtain
\beeqs
\<\vphi_{A_j^\pm}(X,Y),A_{\xi_1}Z\>=0.
\eneqs
Recalling that we have chosen $\xi_1$ such that $\tlamb_i\neq0$ and hence $A_{\xi_1}|_{E_i}$ is an isomorphism of $E_i$, we get that $\vphi_{A_j^\pm}(X,Y)=0$ for all orthonormal pair $X\virg Y\in E_j^\pm$. Lemma \ref{alle} again applies to conclude that $\vphi_{A_j^\pm}$ is also an umbilical bilinear form. Let $Q^\pm\in E_i$ be such that
\beeqs
\vphi_{A_j^\pm}=\<\cdot,\cdot\>Q^\pm.
\eneqs

It remains only to show that $P^\pm$ and $Q^\pm$ vanish. The idea now is to explore how the Codazzi equation relates $P^\pm$ and $Q^\pm$. First, observe that
\beeq\label{suid}
\vphi(A_jX,A_jX)=\sigm_j^2P^-,\hspace{10pt}\vphi(A_jX,X)=Q^-
\eneq
for a unit vector $X\in E_j^+$. One can check these identities by simply writing $X$ as $X=\frac{1}{\sigm_j}A_j^*Y$ with $Y\in E_j^-$ of unit length, since $\frac{1}{\sigm_j}A_j^*:E_j^-\riar E_j^+$ is an orthogonal transformation, and then evaluating the left-hand side using (\ref{AA*s}).

Consider the Codazzi equation for $(X\in E_j^+:\left\|X\right\|=1,A_jX,Z\in E_i^\pm)$,
\beeq\label{lcod}
\al(\lcco_XA_jX,Z)+\al(A_jX,\lcco_XZ)=\al(\lcco_{A_jX}X,Z)+\al(X,\lcco_{A_jX}Z).
\eneq
Taking the inner product of this with $\xi_1$ and using the equation on the right in (\ref{suid}), we have
\beeq\label{cod-lluxipsi}
(\tlamb_j\mp\tlamb_i)\<Q^-,Z\>=-(\tlamb_j\pm\tlamb_i)\<Q^+,Z\>.
\eneq

On the other hand, the Codazzi equation for $(Z\in E_i^\pm,X\in E_j^+:\left\|X\right\|=1,X)$ gives
\beeqs
\noco_Z\al(X,X)-2\al(\lcco_ZX,X)=-\al(\lcco_XZ,X)-\al(Z,\lcco_XX).
\eneqs
Comparing this to the same equation for $(Z\in E_i^\pm,Y\in E_j^-:\left\|Y\right\|=1,Y)$, we see by (\ref{alumeq}) that the two terms involving the normal connection $\noco$ are equal up to sign, so that we can add the equations up in order to get rid of them. After doing so, take the inner product of the resulting equation with $\xi_1$ and $\xi_2$ to obtain
\beeq\label{cod-lluxi1}
(\tlamb_j\pm\tlamb_i)\<P^-,Z\>=(\tlamb_j\mp\tlamb_i)\<P^+,Z\>
\eneq
and
\beeqs
\<P^-,(\rho_j\pm\rho_i)Z+A_i^\pm Z\>-\<P^+,(\rho_j\mp\rho_i)Z-A_i^\pm Z\>=\<Q^++Q^-,Z\>+2\Xi,
\eneqs
respectively, where $\Xi=\<\lcco_ZX,A_jX\>+\<\lcco_ZY,A_j^*Y\>$ is independent of the unit vectors $X\in E_j^+\virg Y\in E_j^-$. In particular, setting $Y=\frac{1}{\sigm_j}A_jX$ (note that $\left\|Y\right\|=1$ by (\ref{AA*s})), we conclude that $\Xi=0$. Therefore,
\beeq\label{cod-llueta2}
\<Q^++Q^-,Z\>=\<P^-,(\rho_j\pm\rho_i)Z+A_i^\pm Z\>-\<P^+,(\rho_j\mp\rho_i)Z-A_i^\pm Z\>.
\eneq
Finally, take the inner product of (\ref{lcod}) with $\xi_2$ and use both equations in (\ref{suid}) to get that
\beeq\label{cod-lluxi2psit}
\<Q^+,(\rho_j\pm\rho_i)Z+A_i^\pm Z\>+\<Q^-,(\rho_j\mp\rho_i)Z-A_i^\pm Z\>=\sigm_j^2\<P^+-P^-,Z\>.
\eneq
Now, if we multiply (\ref{cod-llueta2}) and (\ref{cod-lluxi2psit}) by $(\tlamb_j\pm\tlamb_i)(\tlamb_j\mp\tlamb_i)=(\tlamb_j^2-\tlamb_i^2)$ and use (\ref{cod-lluxipsi}) and (\ref{cod-lluxi1}) into the resulting equations, we obtain a couple of expressions involving only $P^+$ and $Q^+$:
\begin{align}
\label{eq1}\tlamb_i(\tlamb_j\pm\tlamb_i)\<Q^+,Z\>&=\<P^+,(\tlamb_j\mp\tlamb_i)(\rho_j\tlamb_i-\rho_i\tlamb_j)Z\mp\tlamb_j(\tlamb_j\pm\tlamb_i)A_i^\pm Z\>,\\
\label{eq2}\sigm_j^2\tlamb_i(\tlamb_j\mp\tlamb_i)\<P^+,Z\>&=\<Q^+,(\tlamb_j\pm\tlamb_i)(\rho_i\tlamb_j-\rho_j\tlamb_i)Z\pm\tlamb_j(\tlamb_j\mp\tlamb_i)A_i^\pm Z\>.
\end{align}
By changing $Z$ to $A_i^\pm Z$ in (\ref{eq1}) (remind that $A_i^\pm Z\in E_i^\mp$), we obtain
\beeq\label{eq1*}
\tlamb_i(\tlamb_j\mp\tlamb_i)\<Q^+,A_i^\pm Z\>=\<P^+,(\tlamb_j\pm\tlamb_i)(\rho_j\tlamb_i-\rho_i\tlamb_j)A_i^\pm Z\pm\tlamb_j(\tlamb_j\mp\tlamb_i)\sigm_i^2Z\>.
\eneq
Multiplying (\ref{eq2}) by $\tilde{\lamb}_i$ and using (\ref{eq1}) and (\ref{eq1*}) yield an equation just in terms of $P^+$:
\small
\beeq\label{ls1}
((\rho_i\tlamb_j-\rho_j\tlamb_i)^2+\sigm_j^2\tlamb_i^2-\sigm_i^2\tlamb_j^2)(\tlamb_j\mp\tlamb_i)\<P^+,Z\>\pm2(\rho_i\tlamb_j-\rho_j\tlamb_i)\tlamb_j(\tlamb_j\pm\tlamb_i)\<P^+,A_i^\pm Z\>=0.
\eneq
\normalsize
We can again change $Z$ to $A_i^\pm Z$ in (\ref{ls1}), getting
\small
\beeq\label{ls2}
2\sigm_i^2(\rho_i\tlamb_j-\rho_j\tlamb_i)\tlamb_j(\tlamb_j\mp\tlamb_i)\<P^+,Z\>\mp((\rho_i\tlamb_j-\rho_j\tlamb_i)^2+\sigm_j^2\tlamb_i^2-\sigm_i^2\tlamb_j^2)(\tlamb_j\pm\tlamb_i)\<P^+,A_i^\pm Z\>=0.
\eneq
\normalsize
Equations (\ref{ls1}) and (\ref{ls2}) constitute a homogeneous linear system in the variables $\<P^+,Z\>$ and $\<P^+,A_i^\pm Z\>$ whose determinant $d$ is given by
\beeqs
d=\pm(\tlamb_i^2-\tlamb_j^2)[((\rho_i\tlamb_j-\rho_j\tlamb_i)^2+\sigm_j^2\tlamb_i^2-\sigm_i^2\tlamb_j^2)^2+4\sigm_i^2(\rho_i\tlamb_j-\rho_j\tlamb_i)^2\tlamb_j^2].
\eneqs

We show next that $d\neq0$. Suppose that $d=0$. Then, since $\tlamb_j\neq\tlamb_i$,
\begin{align}
	\label{fe1}(\rho_i\tlamb_j-\rho_j\tlamb_i)^2+\sigm_j^2\tlamb_i^2-\sigm_i^2\tlamb_j^2&=0,\\
	\label{fe2}\sigm_i(\rho_i\tlamb_j-\rho_j\tlamb_i)&=0.
\end{align}
Of course, $\sigm_i\neq0$. Otherwise, (\ref{fe1}) and $\tlamb_i\neq0$ would imply that $\sigm_j=0$, which is a contradiction. So, by (\ref{fe2}),
\beeq\label{lfe1}
\rho_i\tlamb_j=\rho_j\tlamb_i.
\eneq
This and (\ref{fe1}) then give
\beeq\label{lfe2}
\sigm_i^2\tlamb_j^2=\sigm_j^2\tlamb_i^2.
\eneq
Now, it follows from Lemma \ref{aebb} and Remark \ref{nola} that $\rho_l^2+\sigm_l^2=\frac{1}{r^2}-\tlamb_l^2$, for $l\in\{i,j\}$. Hence, (\ref{lfe1}) and (\ref{lfe2}) imply that
\beeqs
\big(\frac{1}{r^2}-\tlamb_i^2\big)\tlamb_j^2=\big(\frac{1}{r^2}-\tlamb_j^2\big)\tlamb_i^2,
\eneqs
which leads to a contradiction with $\tlamb_j\neq\tlamb_i$. Therefore, $d\neq0$ and thus $P^+=0$. Finally, (\ref{eq1}) together with (\ref{cod-lluxipsi}) and (\ref{cod-lluxi1}) yields $P^-=0$ and $Q^\pm=0$, as we wished. Hence the lemma is proved.
\enpr

We are now in position to prove Theorem \ref{mine}.

\bepr[Proof of Theorem \ref{mine}]

Let $L$ be a totally geodesic integral submanifold of $E_j$. Since $\al$ is adapted to $(TL,TL^\perp\cap TM)$, it follows from Proposition \ref{roco} that the isometric immersion $f|_L$ admits a reduction of codimension to 2. Moreover, from Lemma \ref{aebb} and (\ref{sofl}) we have that $f|_L$ is minimal. Finally, Remark \ref{xxxx} implies that $f|_L$ also has homothetic Gauss map with the same homothety factor $\frac{1}{r^2}$. Therefore, it follows from Proposition \ref{emes} that $L$ is an Einstein manifold. In other words, $L$ is a minimal Einstein submanifold with codimension two. However, this contradicts Theorem \ref{mats} (resp. Theorem \ref{mihe} if $\dim L=2$), since $f|_L$ is non-totally geodesic. Therefore, $f$ has flat normal bundle and the theorem follows from Nölker's theorem.
\enpr

The following lemma is necessary for the proof of Theorem \ref{mine2}.

\bele\label{app2}
Take an open subset of $M^n$ where $E_1,\dots,E_k$ as in the proof of Theorem \ref{mine} constitute smooth distributions. Then, every $E_i$ such that $\lamb_i(\xi_1)=0$ for some smooth unit normal vector field $\xi_1\in N_fM$ is parallel with respect to the Levi-Civita connection of $M^n$.
\enle
\bepr
Throughout this proof, we take a unit normal vector field $\xi_2$ orthogonal to $\xi_1$ and use the normal frame $\{\xi_1,\xi_2\}$. We consider three cases:

\vspace{\baselineskip}
a) $R^\perp=0$. The assumption that $\lamb_i(\xi_1)=0$ for some smooth unit normal vector field $\xi_1\in N_fM$ is not used in this case. By the Ricci equation, there exists an orthonormal tangent frame $\{X_1,\dots,X_n\}$ such that
\beeqs
\al(X_i,X_j)=0\virg1\leq i\neq j\leq n.
\eneqs
Therefore, for each $x\in M^n$ the tangent space $T_xM$ decomposes orthogonally as
\beeqs
T_xM=D_1(x)\oplus\dots\oplus D_s(x),
\eneqs
where each $D_i(x)$ is a common eigenspace of all shape operators, that is,
\beeqs
A_\xi X_i=\mu_i(\xi)X_i
\eneqs
if $X_i\in D_i(x)\virg1\leq i\leq s=s(x)$, and $\mu_i\neq\mu_j$ for $i\neq j$. Now, it follows by the uniqueness part of Lemma \ref{aldec} that $s=k\virg \mu_i=\lamb_i$ and $D_i=E_i\virg1\leq i\leq k$. In this special case, the maps $\xi\mapsto\lamb_i(\xi)$ are linear and hence there exist unique normal vector fields $\eta_i\virg1\leq i\leq k$, called the \emph{principal normals} of $f$, such that $\lamb_i(\xi)=\<\eta_i,\xi\>$. Therefore,
\beeqs
E_i=\{X\in TM:\al(X,Y)=\<X,Y\>\eta_i\mbox{ for all }Y\in TM\}\virg 1\leq i\leq k,
\eneqs
and the second fundamental form of $f$ has the simple representation
\beeq\label{app2sr}
\al(X,Y)=\sum_{i=1}^{k}\<X^i,Y^i\>\eta_i
\eneq
for all $X\virg Y\in TM$, where $X\mapsto X^i$ is the orthogonal projection onto $E_i$. Then, the assumption on the Gauss map implies that
\beeq\label{app2a}
\left\|\eta_i\right\|^2=\thff(X_i,X_i)=\frac{1}{r^2},
\eneq
where $X_i\in E_i$ is a unit vector. Therefore, since $\eta_i\neq\eta_j\virg1\leq i\neq j\leq k$, it follows from the Cauchy-Schwarz inequality that
\beeq\label{app2acs}
\<\eta_i,\eta_j\><\frac{1}{r^2}\virg1\leq i\neq j\leq k.
\eneq

Consider the tensor $\varphi_{ij}:TM\times E_i\riar E_j$ defined by $\varphi_{ij}(X,Y)=(\lcco_XY)_{E_j}\virg1\leq i\neq j\leq k$. To conclude that $E_i$ is parallel in the Levi-Civita connection, $1\leq i\leq k$, we must show that all $\varphi_{ij}$ are identically zero, for $1\leq i\neq j\leq k$.

The Codazzi equation for $(Z\in E_j,X\in E_i, Y\in E_i)$ and (\ref{app2sr}) give
\beeq\label{appbnoco}
\<X,Y\>\noco_Z\eta_i=\<\varphi_{ij}(X,Y),Z\>(\eta_i-\eta_j).
\eneq
Taking the inner product with $\eta_i$, we have, by (\ref{app2a}) and (\ref{app2acs}),
\beeq\label{app2atg}
\<\varphi_{ij}(X,Y),Z\>=0.
\eneq
Since $X\virg Y\in E_i\virg Z\in E_j$ and the indices $i\neq j$ have been arbitrarily chosen, the above equation implies that each $E_i$ is a totally geodesic distribution, $1\leq i\leq k$. Thus, in order to conclude that $E_i$ is parallel in the Levi-Civita connection, it remains only to check (\ref{app2atg}) for $X\in E_l\virg Y\in E_i\virg Z\in E_j$ and pairwise distinct indices $i\virg j\virg l$, since $\<\vphi_{ij}(X,Y),Z\>=-\<\vphi_{ji}(X,Z),Y\>=0$ for $X\virg Z\in E_j\virg Y\in E_i$.

We claim that this follows from the Codazzi equation for $(X\in E_l,Y\in E_i,Z\in E_j)$. In fact, the latter gives
\beeqs
\<\varphi_{ij}(X,Y),Z\>(\eta_j-\eta_i)=\<\varphi_{lj}(Y,X),Z\>(\eta_j-\eta_l).
\eneqs
But since $\eta_i\virg\eta_j\virg\eta_l$ are pairwise distinct and have the same norm, it is straightforward to conclude that the vectors $\eta_j-\eta_i$ and $\eta_j-\eta_l$ cannot be collinear, so that 
\beeq\label{appbij}
\<\varphi_{ij}(X,Y),Z\>=0,
\eneq
as we wished.

\vspace{\baselineskip}
b) $\dim E_i\geq2$. We show that $R^\perp=0$, hence reducing the problem to the previous case. Since $\lamb_i(\xi_1)=0$, we have that $E_i\subseteq\ker A_{\xi_1}$. Furthermore, the assumption of homothetic Gauss map implies that
\beeq\label{app2hgm}
A_{\xi_2}^2|_{\ker A_{\xi_1}}=\frac{1}{r^2}I_{\ker A_{\xi_1}}.
\eneq
From this we then obtain that
\beeqs
A_\xi^2|_{\ker A_{\xi_1}}=\big(\frac{\<\xi,\xi_2\>}{r}\big)^2I_{\ker A_{\xi_1}}
\eneqs
for every $\xi\in N_fM$, so that $\ker A_{\xi_1}$ fits into decomposition (\ref{alde}). By uniqueness, we conclude that actually $E_i=\ker A_{\xi_1}$. Now, it is a consequence of (\ref{A2}) and (\ref{app2hgm}) that
\beeqs
E_i\subseteq E_{A_{\xi_2}}\big(-\frac{1}{r}\big)\oplus E_{A_{\xi_2}}\big(\frac{1}{r}\big).
\eneqs
We claim that equality holds in the above inclusion. Indeed, take for instance a vector $X\in E_{A_{\xi_2}}\big(\frac{1}{r}\big)$. In particular, $A_{\xi_2}^2X=\frac{1}{r^2}X$. The assumption on the Gauss map then yields $A_{\xi_1}^2X=0$
and, consequently, $X\in \ker A_{\xi_1}=E_i$. In other words, $E_{A_{\xi_2}}\big(\frac{1}{r}\big)\subset E_i$. Similarly, we show that $E_{A_{\xi_2}}\big(-\frac{1}{r}\big)\subset E_i$, so that our claim is proved, i.e.,
\beeq\label{app2es}
E_i=E_{A_{\xi_2}}\big(-\frac{1}{r}\big)\oplus E_{A_{\xi_2}}\big(\frac{1}{r}\big).
\eneq

Take an orthonormal frame $\{X_1,\dots,X_m\}$ of eigenvectors of $A_{\xi_2}|_{E_i}$, so that $A_{\xi_2}X_j=\pm\frac{1}{r}X_j$, $1\leq j\leq m$. Note that $\al(X_j,X_l)=\pm\frac{1}{r}\delta_{jl}\xi_2$. Let $\omega$ be the normal connection 1-form on $TM$ defined by $\omega(X)=\<\noco_X\xi_1,\xi_2\>$. We will check that $\omega=0$ to conclude that $R^\perp=0$, since the codimension is two.

We claim that $E_i$ is a totally geodesic distribution. To see this, consider the tensor $\varphi:E_i\times E_i\riar E_i^\perp\cap TM$ defined by $\varphi(X,Y)=(\lcco_XY)_{E_i^\perp\cap TM}$. It suffices to show that $\varphi(X_j,X_l)=0$ for $1\leq j\virg l\leq m$. The Codazzi equation for $(Y\in E_i^\perp\cap TM,X_j,X_j)$ yields
\beeq\label{app2ce}
\pm\frac{1}{r}\noco_Y\xi_2=-\al(\lcco_{X_j}Y,X_j)-\al(Y,\lcco_{X_j}X_j).
\eneq
Taking the inner product with $\xi_2$ gives $\big\<\varphi(X_j,X_j),\big(A_{\xi_2}\mp\frac{1}{r}I_{TM}\big)Y\big\>=0$. But, by (\ref{app2es}), we have that $A_{\xi_2}\mp\frac{1}{r}I_{TM}$ maps $E_i^\perp\cap TM$ onto $E_i^\perp\cap TM$. Thus, we conclude from the above that $\varphi(X_j,X_j)=0$
for all $1\leq j\leq m$. Take now the inner product of (\ref{app2ce}) with $\xi_1$. By the above and $E_i=\ker A_{\xi_1}$ we obtain
\beeq\label{app2nc1}
\omega(Y)=0
\eneq
for all $Y\in E_i^\perp\cap TM$.

On the other hand, the Codazzi equation for $(Y\in E_i^\perp\cap TM,X_j,X_l)\virg j\neq l$, now gives
\beeqs
\al(\lcco_YX_j,X_l)+\al(X_j,\lcco_YX_l)=\al(\lcco_{X_j}Y,X_l)+\al(Y,\lcco_{X_j}X_l).
\eneqs
Since $E_i=\ker A_{\xi_1}$, taking the inner product with $\xi_1$ yields
\beeqs
\<\varphi(X_j,X_l),A_{\xi_1}Y\>=0.
\eneqs
However, $A_{\xi_1}|_{E_i^\perp\cap TM}$ is an isomorphism of $E_i^\perp\cap TM$, and therefore
\beeqs
\varphi(X_j,X_l)=0.
\eneqs
for $j\neq l$. This concludes the proof of the claim.

Finally, Codazzi equation for $(X\in E_i, Y\in E_i, \xi_1)$ together with the claim just proved implies that
\beeqs
\omega(Y)A_{\xi_2}X=\omega(X)A_{\xi_2}Y.
\eneqs
Since $A_{\xi_2}|_{E_i}:E_i\riar E_i$ is an isomorphism and we are under the assumption $\dim E_i\geq2$, it follows that
\beeqs
\omega|_{E_i}=0.
\eneqs
This and (\ref{app2nc1}) show that $\omega$ vanishes identically and thus $R^\perp=0$, as we wished.

\vspace{\baselineskip}
c) Neither a) nor b) occurs. Let $\varphi_{jl}:TM\times E_j\riar E_l$ be the tensor defined as in case a), for any pair of distinct indices $j\virg l$. Set
\beeqs
\Gamma=\{j:\lamb_j(\xi)=0\mbox{ for some }\xi\in N_fM\}.
\eneqs
By assumption, $i\in\Gamma$. If $\dim E_j\geq2$ for some $j\in\Gamma$, we can conclude as in case b) that $R^\perp=0$. Therefore, there is no loss of generality in assuming that all $E_j$ for $j\in\Gamma$ are line bundles. Let $E_j$ be locally spanned by a unit vector field $X_j$. So, $\{X_j:j\in\Gamma\}$ is an orthonormal basis of $F=\oplus_{j\in\Gamma}E_j$ that diagonalizes all shape operators. Then, we can use the same argument as in case a) to show that, if $j\in\Gamma$,
\beeq\label{app2me1}
\varphi_{ij}(X,Y)=0
\eneq
for all $X\in F$ and $Y\in E_i$. To check that the same holds for $X\notin F$, we can assume by tensoriality that $X\in E_l$ with $l\notin\Gamma$, so that the second fundamental form restricted to $E_l$ has the algebraic structure given by Lemma \ref{aebb}. For simplicity, we set $\lamb_j(\xi_1)=\tlamb_j$ and $\lamb_j(\xi_2)=\rho_j$ for $j\in\Gamma$. Notice that $\tlamb_i=0$ by assumption and $\rho_i=\pm\frac{1}{r}$ by Remark \ref{nola}. Replacing $\xi_2$ by $-\xi_2$ if necessary, we can assume $\rho_i=\frac{1}{r}$. Furthermore, it holds that $\tlamb_j\neq0$, for $E_i=\ker A_{\xi_1}$.

Using the Codazzi equation for $(X_i,X_j,X\in E_l^\pm )$ we have
\beeqs
\al(\lcco_{X_i}X_j,X)+\al(X_j,\lcco_{X_i}X)=\al(\lcco_{X_j}X_i,X)+\al(X_i,\lcco_{X_j}X).
\eneqs
Taking the inner product with $\xi_1$ and $\xi_2$ yields
\beeq\label{app2cijlxi1}
(\tlamb_j\mp\tlamb_l)\<\varphi_{jl}(X_i,X_j),X\>=\mp\tlamb_l\<\varphi_{il}(X_j,X_i),X\>
\eneq
and
\beeq\label{app2cijlxi2}
\<\varphi_{jl}(X_i,X_j),(A_l^\pm -(\rho_j\mp\rho_l)I_{E_l})X\>=\big\<\varphi_{il}(X_j,X_i),\big(A_l^\pm -\big(\frac{1}{r}\mp\rho_l\big)I_{E_l}\big)X\big\>,
\eneq
respectively. Multiplying (\ref{app2cijlxi2}) by $\tlamb_l$ and using (\ref{app2cijlxi1}), we obtain
\beeq\label{app2cijlf}
\big\<\varphi_{jl}(X_i,X_j),\big(\tlamb_j A_l^\pm +\big(\tlamb_j\big(\frac{1}{r}\mp\rho_l\big)\mp\tlamb_l\big(\frac{1}{r}-\rho_j\big)\big)I_{E_l}\big)X\big\>=0.
\eneq

Now, using the Codazzi equation for $(X\in E_l^\pm ,X_i,X_j)$ we have
\beeqs
\al(\lcco_XX_i,X_j)+\al(X_i,\lcco_XX_j)=\al(\lcco_{X_i}X,X_j)+\al(X,\lcco_{X_i}X_j).
\eneqs
On the other hand, taking the inner product with $\xi_1$ and $\xi_2$ we get
\beeq\label{app2clijxi1}
(\tlamb_j\mp\tlamb_l)\<\varphi_{jl}(X_i,X_j),X\>=-\tlamb_j\<\varphi_{ij}(X,X_i),X_j\>
\eneq
and
\beeq\label{app2clijxi2}
\big(\frac{1}{r}-\rho_j\big)\<\varphi_{ij}(X,X_i),X_j\>=\<\varphi_{jl}(X_i,X_j),((\rho_j\mp\rho_l)I_{E_l}-A_l^\pm )X\>,
\eneq
respectively, where we set for convenience $\varphi_{ij}(X,X_i)=0$ in the case $i=j$. Multiplying (\ref{app2clijxi2}) by $\tlamb_j$ and using (\ref{app2clijxi1}) give
\beeq\label{app2clijf}
\big\<\varphi_{jl}(X_i,X_j),\big(\tlamb_j A_l^\pm -\big(\tlamb_j\big(\frac{1}{r}\mp\rho_l\big)\mp\tlamb_l\big(\frac{1}{r}-\rho_j\big)\big)I_{E_l}\big)X\big\>=0.
\eneq
Finally, add (\ref{app2cijlf}) and (\ref{app2clijf}) to conclude that $\<\varphi_{jl}(X_i,X_j),A_l^\pm X\>=0$ for all $X\in E_l^\pm $. Since $A_l^\pm :E_l^\pm \riar E_l^\mp $ is an isomorphism, it follows that
\beeq\label{app2idk}
\varphi_{jl}(X_i,X_j)=0.
\eneq
This together with (\ref{app2clijxi1}) shows that (\ref{app2me1}) also holds for $X\in E_l^\pm \virg l\notin\Gamma$, and hence $\varphi_{ij}=0$ for every $j\in\Gamma$. It remains to verify that $\varphi_{il}=0$ for $l\notin\Gamma$. But then we know from Lemma \ref{Ejtg} that $E_l$ is a totally geodesic distribution. In particular,
\beeq\label{app2iat}
\varphi_{il}(X,X_i)=0\virg\forall X\in E_l.
\eneq
Moreover, it follows from (\ref{app2cijlxi1}) and (\ref{app2idk}) that (\ref{app2iat}) also holds for $X=X_j$ with $j\in\Gamma$, and thus for all $X\in F$. So, in order to conclude that $\varphi_{il}=0$, it remains only to check (\ref{app2iat}) for $X\in E_{l'}$ with $l'\notin\Gamma$ and $l'\neq l$.

From the Codazzi equation for $(X\in E_{l'}^+,X_i,Y\in E_l^+)$, we obtain
\beeqs
\al(\lcco_XX_i,Y)+\al(X_i,\lcco_XY)=\al(\lcco_{X_i}X,Y)+\al(X,\lcco_{X_i}Y).
\eneqs
Taking the inner product with $\xi_1$ and $\xi_2$ yields
\beeq\label{app2lijnxi1}
\tlamb_l\<\varphi_{il}(X,X_i),Y\>=(\tlamb_l-\tlamb_{l'})\<\varphi_{l'l}(X_i,X),Y\>
\eneq
and
\beeqs
\begin{aligned}
\big\<\varphi_{il}(X,X_i),\big(A_l-\big(\frac{1}{r}-\rho_l\big)I_{E_l^+}\big)Y\big\>=\<\varphi_{l'l}(X_i,X),(A_l-(\rho_{l'}-\rho_l)I_{E_l^+})Y & \> \\
         +\<\varphi_{ll'}(X_i,Y),A_{l'} X & \>,
\end{aligned}
\eneqs
respectively. A similar computation for $X\in E_{l'}^{\epsilon'}\virg Y\in E_l^{\epsilon}$, where $\epsilon\virg\epsilon'\in\{+,-\}$, gives
\beeq\label{app2lijnff}
\epsilon'\tlamb_{l'}\<\varphi_{l'l}(X_i,X),A_l^{\epsilon} Y\>=\epsilon\tlamb_l\<\varphi_{ll'}(X_i,Y),A_{l'}^{\epsilon'}X\>.
\eneq
Now, multiply the above equation by $\tlamb_l$ and use (\ref{app2lijnff}), to get
\beeqs
\begin{aligned}
\big(\tlamb_l\big(\frac{1}{r}-\rho_{l'}\big)-\tlamb_{l'}\big(\frac{1}{r}-\rho_l\big)\big)\<\varphi_{l'l}(X_i,X),Y\>=\tlamb_{l'}\<\varphi_{l'l}(X_i,X),A_l Y & \> \\
    -\tlamb_l\<\varphi_{ll'}(X_i,Y),A_{l'} X & \>.
\end{aligned}
\eneqs
If we then invert the roles of $l$ and $l'$ (and $X$ and $Y$) in the above equation, we see that, while the left-hand side remains the same, since $\<\vphi_{l'l}(X_i,X),Y\>=-\<\vphi_{ll'}(X_i,Y),X\>$, the right-hand side changes sign. Therefore, both sides must vanish, i.e.,
\beeqs
\big(\tlamb_l\big(\frac{1}{r}-\rho_{l'}\big)-\tlamb_{l'}\big(\frac{1}{r}-\rho_l\big)\big)\<\varphi_{l'l}(X_i,X),Y\>=0
\eneqs
and
\beeq\label{app2lijnf}
\tlamb_{l'}\<\varphi_{l'l}(X_i,X),A_lY\>=\tlamb_l\<\varphi_{ll'}(X_i,Y),A_{l'} X\>.
\eneq
Suppose, by contradiction, that $\<\varphi_{l'l}(X_i,X),Y\>\neq0$ for certain $X\in E_{l'}^+$ and $Y\in E_l^+$. Then
\beeq\label{app2lijnfe}
\tlamb_l\big(\frac{1}{r}-\rho_{l'}\big)=\tlamb_{l'}\big(\frac{1}{r}-\rho_l\big).
\eneq
Set $X=A_{l'}\tilde{X}\virg Y=A_l^*\tilde{Y}$ in (\ref{app2lijnff}) and recall (\ref{AA*s}), obtaining
\beeqs
\tlamb_{l'}\sigm_l^2\<\varphi_{ll'}(X_i,\tilde{Y}),A_{l'}\tilde{X}\>=-\tlamb_l\sigm_{l'}^2\<\varphi_{l'l}(X_i,\tilde{X}),A_l^*\tilde{Y}\>.
\eneqs
This and (\ref{app2lijnff}) give
\beeqs
((\tlamb_{l'}\sigm_l)^2-(\tlamb_l\sigm_{l'})^2)\<\varphi_{l'l}(X_i,\tilde{X}),A_l^*\tilde{Y}\>=0.
\eneqs
However, since we are under the assumption that $\<\varphi_{l'l}(X_i,X),Y\>\neq0$ for certain $X\in E_{l'}^+\virg Y\in E_l^+$ and $A_l^*$ is onto $E_l^+$, it follows that
\beeq\label{happynewyear}
(\tlamb_l\sigm_{l'})^2=(\tlamb_{l'}\sigm_l)^2.
\eneq
This together with $\tlamb_l^2+\rho_l^2+\sigm_l^2=\frac{1}{r^2}$ (and the same for $l'$) implies that
\beeqs
\tlamb_l^2\big(\frac{1}{r^2}-\rho_{l'}^2\big)=\tlamb_{l'}^2\big(\frac{1}{r^2}-\rho_l^2\big).
\eneqs
This and (\ref{app2lijnfe}) imply
\beeq\label{happy2014}
\tlamb_l^2\big(\frac{1}{r}-\rho_{l'}\big)=\tlamb_{l'}^2\big(\frac{1}{r}-\rho_l\big).
\eneq
Comparing it to (\ref{app2lijnfe}), we finally obtain $\tlamb_l=\tlamb_{l'}$. But then (\ref{happynewyear}) and (\ref{happy2014}) yield $\rho_l=\rho_{l'}$ and $\sigm_l=\sigm_{l'}$. However, those three relations imply that $E=E_l\oplus E_{l'}$ fits into decomposition (\ref{alde}), which contradicts its uniqueness. Therefore, we have that $\<\varphi_{l'l}(X_i,X),Y\>=0$ for all $X\in E_{l'}^+\virg Y\in E_l^+$. Finally, (\ref{app2lijnxi1}) then implies that $\<\varphi_{il}(X,X_i),Y\>=0$ for every $X\in E_{l'}^+\virg Y\in E_l^+$. Entirely analogous arguments give $\<\varphi_{il}(X,X_i),Y\>=0$ for all $X\in E_{l'}^\pm $, $Y\in E_l^\pm $, and therefore $\varphi_{il}=0$. This completes the proof of Lemma \ref{app2}.
\enpr

Observe that the proofs of Lemmas \ref{Ejtg} and \ref{app2} make only use of the Codazzi equation, which is the same for space forms of nonzero curvature. Thus, we conclude that the lemmas remain true in this setting. This will be used in the proof of Theorem \ref{mine2}.

\bepr[Proof of Theorem \ref{mine2}]
For the hyperbolic space, the previous proof works \emph{mutatis mutandis}, since Theorem \ref{mihe} and Theorem \ref{mats} are also true in this setting. The situation for the sphere is more delicate. By Lemma \ref{app2}, every distribution $E_i$ such that $\lamb_i(\xi)=0$ for some smooth unit normal vector field $\xi\in N_fM$ is parallel with respect to the Levi-Civita connection of $M^n$. But since $\al$ is adapted to (\ref{alded}) and we are under the assumption that $f$ is substantial and irreducible, it follows that no such $E_i$ must appear. So, we conclude that the only blocks $E_i$ composing (\ref{alded}) are those to which Lemma \ref{aebb} applies. They are all `minimal blocks' in the sense that the trace of any $A_\xi$ restricted to $E_i$ is zero. Therefore, $f$ itself must be a minimal isometric immersion, and consequently $M^n$ is Einstein by Proposition \ref{emes}, which is valid in any space form. But since all the possibilities in Matsuyama's classification are reducible, then $n=2$ and $f(M^n)$ is a piece of the Veronese surface, by Kenmotsu's result \cite{MR706526}, since the Clifford torus is also reducible. 
\enpr

\bere\label{reapp2}
Case a) in the proof of Lemma \ref{app2} is the quintessence of Nölker's argument to prove his theorem. Indeed, first observe that the proof works for arbitrary codimension. Then by the fact that every $E_i\virg1\leq i\leq k$, is a parallel distribution and de Rham's theorem, $M^n$ is the Riemannian product of the integral manifolds $M_1,\dots,M_k$ of $E_1,\dots,E_k$ through one point $p_0\in M^n$. Since $\al$ is adapted to the product net $(E_1,\dots,E_k)$, $f$ is a Riemannian product of isometric immersions $f_i:M_i\riar\R^{n_i}\virg i=1,\dots,k$, by the well-known lemma of Moore (\cite{MR0307128}, p. 163). From (\ref{app2sr}) and (\ref{app2a}) follows that $f_1,\dots,f_k$ are totally umbilical immersions with mean curvature vectors of constant length, thus Euclidean spheres or curves of constant curvature.
\enre

In light of the results presented so far, we conclude this section posing the following conjecture suggesting a possible complete solution to our Main Problem in arbitrary codimension.

\beco
Let $f:M^n\riar\Q_c^{n+p}$ be an irreducible isometric immersion with homothetic Gauss map, $n\geq2$. Then $M^n$ is an Einstein manifold and, up to composition with a totally umbilical inclusion, $f$ is a minimal isometric immersion into some sphere $\Sp_{\tilde{c}}^{n+q}\virg q\leq p$.
\enco

\bere
The preceding conjecture is stronger than Conjecture \ref{disc} and also implies the version of the latter for hyperbolic space forms. The conjecture is true for compact orientable Einstein submanifolds of Euclidean space whose Gauss map is harmonic, according to a result due to Mut$\bar{\mbox{o}}$ \cite{MR575405}.
\enre

\bibliographystyle{plain}

\bibliography{mybib2}

\begin{thebibliography}{10}

\bibitem{MR787964}
Robert~L. Bryant.
\newblock Minimal surfaces of constant curvature in {$S^n$}.
\newblock {\em Trans. Amer. Math. Soc.}, 290(1):259--271, 1985.

\bibitem{MR0057000}
Eugenio Calabi.
\newblock Isometric imbedding of complex manifolds.
\newblock {\em Ann. of Math. (2)}, 58:1--23, 1953.

\bibitem{MR749575}
Bang-Yen Chen.
\newblock {\em Total mean curvature and submanifolds of finite type}, volume~1
  of {\em Series in Pure Mathematics}.
\newblock World Scientific Publishing Co., Singapore, 1984.

\bibitem{MR0226514}
Shiing-shen Chern and Robert Osserman.
\newblock Complete minimal surfaces in euclidean {$n$}-space.
\newblock {\em J. Analyse Math.}, 19:15--34, 1967.

\bibitem{MR1075013}
Marcos Dajczer.
\newblock {\em Submanifolds and isometric immersions}, volume~13 of {\em
  Mathematics Lecture Series}.
\newblock Publish or Perish Inc., Houston, TX, 1990.
\newblock Based on the notes prepared by Mauricio Antonucci, Gilvan Oliveira,
  Paulo Lima-Filho and Rui Tojeiro.

\bibitem{MR2000031}
Antonio~J. Di~Scala.
\newblock Minimal immersions of {K}\"ahler manifolds into {E}uclidean spaces.
\newblock {\em Bull. London Math. Soc.}, 35(6):825--827, 2003.

\bibitem{MR0278318}
Manfredo~P. do~Carmo and Nolan~R. Wallach.
\newblock Minimal immersions of spheres into spheres.
\newblock {\em Ann. of Math. (2)}, 93:43--62, 1971.

\bibitem{MR706526}
Katsuei Kenmotsu.
\newblock Minimal surfaces with constant curvature in {$4$}-dimensional space
  forms.
\newblock {\em Proc. Amer. Math. Soc.}, 89(1):133--138, 1983.

\bibitem{MR1255769}
Yoshio Matsuyama.
\newblock Minimal {E}instein submanifolds with codimension two.
\newblock {\em Tensor (N.S.)}, 52(1):61--68, 1993.

\bibitem{MR0307128}
John~Douglas Moore.
\newblock Isometric immersions of riemannian products.
\newblock {\em J. Differential Geometry}, 5:159--168, 1971.

\bibitem{MR575405}
Yosio Muto.
\newblock Submanifolds of a {E}uclidean space with homothetic {G}auss map.
\newblock {\em J. Math. Soc. Japan}, 32(3):531--555, 1980.

\bibitem{MR1066578}
Stefan N{\"o}lker.
\newblock Isometric immersions with homothetical {G}auss map.
\newblock {\em Geom. Dedicata}, 34(3):271--280, 1990.

\bibitem{MR0234388}
Morio Obata.
\newblock The {G}auss map of immersions of {R}iemannian manifolds in spaces of
  constant curvature.
\newblock {\em J. Differential Geometry}, 2:217--223, 1968.

\bibitem{MR775143}
Yoshihiro Ohnita.
\newblock The first standard minimal immersions of compact irreducible
  symmetric spaces.
\newblock In {\em Differential geometry of submanifolds ({K}yoto, 1984)},
  volume 1090 of {\em Lecture Notes in Math.}, pages 37--49. Springer, Berlin,
  1984.

\bibitem{MR609562}
Robert Osserman.
\newblock Minimal surfaces, {G}auss maps, total curvature, eigenvalue
  estimates, and stability.
\newblock In {\em The {C}hern {S}ymposium 1979 ({P}roc. {I}nternat. {S}ympos.,
  {B}erkeley, {C}alif., 1979)}, pages 199--227. Springer, New York, 1980.

\bibitem{MR0407774}
Nolan~R. Wallach.
\newblock Minimal immersions of symmetric spaces into spheres.
\newblock In {\em Symmetric spaces ({S}hort {C}ourses, {W}ashington {U}niv.,
  {S}t. {L}ouis, {M}o., 1969--1970)}, pages 1--40. Pure and Appl. Math., Vol.
  8. Dekker, New York, 1972.

\end{thebibliography}

\vspace{\baselineskip}
Instituto Nacional de Matemática Pura e Aplicada (IMPA), Estrada Dona Castorina 110, 22460-320 Rio de Janeiro, Brazil; gfreitas@impa.br

\end{document}